\documentclass[11pt]{article}
\usepackage[margin=1.05in]{geometry}
\usepackage[T1]{fontenc}
\usepackage{lmodern}
\usepackage{amsmath,amssymb,amsthm,mathtools}
\usepackage{microtype,booktabs,array,enumitem}
\usepackage[colorlinks=true,linkcolor=blue!45!black,citecolor=blue!45!black,urlcolor=blue!45!black]{hyperref}
\usepackage{xcolor}
\hypersetup{pdftitle={The Lorentz Threshold in Fourier Restriction},pdfauthor={Rishad Shahmurov and Veli Shahmurov},pdfkeywords={Stein--Tomas restriction; Gaussian type; Lorentz coefficient spaces; capillary waves; gamma-radonifying operators}}
\newtheorem{theorem}{Theorem}[section]
\newtheorem{lemma}[theorem]{Lemma}
\newtheorem{proposition}[theorem]{Proposition}
\newtheorem{corollary}[theorem]{Corollary}
\newtheorem{remark}[theorem]{Remark}
\newtheorem{conjecture}[theorem]{Conjecture}
\newtheorem{problem}[theorem]{Problem}
\newcommand{\R}{\mathbb R}

\newcommand{\E}{\mathbb E}
\newcommand{\ga}{\gamma}
\newcommand{\calS}{\mathsf S}
\newcommand{\calF}{\mathsf F}
\newcommand{\one}{\mathbf 1}
\newcommand{\norm}[1]{\left\|#1\right\|}
\title{\textbf{The Lorentz Threshold in Fourier Restriction}}
\author{Rishad Shahmurov\thanks{Cellular Products Research and Development, 595 East Crossville Rd., Roswell, GA 30075, USA. Correspondence: \texttt{rshahmurov@crimson.ua.edu}.}
\and Veli Shahmurov\thanks{Antalya Bilim University, Department of Industrial Engineering, 07190 Antalya, Turkey; Azerbaijan State Economic University, AZ1001 Baku, Azerbaijan. Emails: \texttt{veli.sahmurov@antalya.edu.tr}, \texttt{veli.sahmurov@gmail.com}.}}
\date{}
\begin{document}
\maketitle

\begin{abstract}
Classical Stein--Tomas restriction controls a Fourier transform on a curved
surface, but the estimate changes when its values lie in a Banach space.
We determine the sharp Lorentz threshold at the critical exponent $4/3$:
for coefficients in $\ell_D^{4/3,q}$, the optimal Gaussian restriction
constant grows like $(1+\log D)^{|3/4-1/q|}$. A multiscale wave-packet
construction proves sharpness, including the delicate index $q=2$.
The same threshold holds for fractional Gaussian noise, while Gaussian
type and cotype alone do not distinguish the positive and negative cases.
As a dispersive consequence, we establish a sharp Lorentz-channel
spacetime estimate for linear capillary waves.
\end{abstract}


\section{Introduction}

\subsection{Why restricting a Fourier transform is difficult}

The Fourier transform describes how a function is assembled from different
frequencies. Restricting it to a surface sounds like ordinary evaluation, but
there is a basic difficulty: a sphere has zero three-dimensional volume, and
the Fourier transform of a general $L^p$ function is not defined
pointwise on that sphere. Fourier restriction theory asks when the shape of
the surface makes such a trace meaningful.

A fundamental answer is the Stein--Tomas theorem. For a compact curved
hypersurface $\Sigma\subset\R^3$, it gives
\begin{equation}\label{eq:scalar-st}
 \|\widehat f|_\Sigma\|_{L^2(\Sigma)}
 \le C_\Sigma \|f\|_{L^{4/3}(\R^3)}.
\end{equation}
Curvature is essential: oscillations from different points of the surface
partly cancel, allowing a Fourier trace that would not follow from a simple
pointwise estimate. The restriction exponent $4/3$ is critical. A narrow
frequency cap corresponds to a long, thin physical wave packet. For the
packet $g_R$ defined in \eqref{eq:packet} below,
\[
 \|g_R\|_{L^s(\R^3)}\asymp R^{4/s-3},
 \qquad \|\widehat g_R|_\Sigma\|_2\gtrsim1.
\]
Thus $s>4/3$ cannot replace $4/3$ uniformly. This elementary scaling
already explains why the endpoint is particularly delicate.

The estimate emerged from the classical work of Stein and Tomas
\cite{Tomas,SteinBook}; its importance extends to Fourier extension,
oscillatory integrals and dispersive-wave estimates. These connections and
more recent endpoint developments are reviewed in \cite{OliveiraSurvey}.
Bak and Seeger \cite{BakSeeger} sharpened the scalar estimate at the same
critical exponent by replacing its physical source with a Lorentz space,
$L^{4/3,2}(\R^3)$. Our question differs: the physical exponent remains
$4/3$, while the \emph{values} of the functions belong to a Lorentz space.

\subsection{Many coefficient channels and Gaussian type}

Imagine that $f(x)$ records not one scalar measurement but a vector of many
channels. Applying \eqref{eq:scalar-st} to each channel separately does not
answer whether the joint estimate stays stable as the number of channels
grows. The answer depends on the Banach space $X$ in which the values live.

The natural joint output is a $\gamma$-radonifying operator. For a finite-rank
operator $T:H\to X$ and an orthonormal system $(h_\nu)$ spanning its Hilbert
input, its norm is
\begin{equation}\label{eq:intro-gamma}
 \|T\|_{\gamma(H,X)}
 =\left(\E\left\|\sum_\nu\gamma_\nu Th_\nu\right\|_X^2\right)^{1/2},
\end{equation}
where the $\gamma_\nu$ are independent standard Gaussians. Informally,
\eqref{eq:intro-gamma} measures the size of all channels after Gaussian
mixing. When $X$ is Hilbert, it reduces to the usual square-sum of channel
energies; for a general Banach space it retains finer geometric information.

Fix a compact $C^2$ surface $\Sigma$ with nonvanishing Gaussian curvature
and a smooth positive density $d\sigma$. We study
\begin{equation}\label{eq:intro-frt}
 R_\Sigma:L^{p}(\R^3;X)\longrightarrow
 \gamma(L^2(\Sigma),X),\qquad p=\frac43,
\end{equation}
and denote its norm by $\calS_p(X;\Sigma)$. For finite-rank Schwartz
functions, $R_\Sigma f$ is simply the Fourier trace $\widehat f|_\Sigma$
interpreted as an operator on $L^2(\Sigma)$.

A natural first prediction comes from \emph{Gaussian type}. A Banach space
has Gaussian type $p$ if
\begin{equation}\label{eq:intro-type}
 T_p(X):=\sup_{(x_j)\ne0}
 \frac{\left(\E\|\sum_j\gamma_jx_j\|_X^2\right)^{1/2}}
 {\left(\sum_j\|x_j\|_X^p\right)^{1/p}}<\infty.
\end{equation}
Type $p$ says that independent Gaussian contributions cannot accumulate
faster than their $\ell^p$ input sizes. The endpoint restriction estimate
forces type $p$, but its converse is much less clear: the Fourier traces of
many physical packets are \emph{correlated}, not independent. Work on
Banach-valued Gaussian embeddings \cite{KNVW,Veraar} explains why type is a
natural candidate, while leaving the curved endpoint as an additional
problem. Other recent vector-valued restriction results can be found in
\cite{SchattenCompanion}.

\subsection{The sharp Lorentz threshold}

Lorentz spaces provide an unusually precise test of this question. Let
$1\le q<\infty$, $D\ge2$, and put $E_D^{p,q}=\ell_D^{p,q}$, using a Banach
lattice norm equivalent to
\begin{equation}\label{eq:lorentz}
 \|z\|_{p,q}^*=
 \left(\sum_{r=1}^{D}r^{q/p-1}(z^*(r))^q\right)^{1/q},
\end{equation}
where $z^*$ denotes the decreasing rearrangement of the absolute values.
The parameter $p$ controls the main integrability scale, while $q$ measures
how values from different levels combine. For $q>p$ the displayed
expression may be a quasi-norm; the equivalent Banach norm is used
throughout, with constants independent of $D$.

Our first theorem gives an exact quantitative answer.

\begin{theorem}[Sharp Lorentz restriction threshold]\label{thm:main}
Let $p=4/3$, $1\le q<\infty$, and $\Sigma$ be as above. There exist
constants $0<c_{\Sigma,q}\le C_{\Sigma,q}<\infty$ such that for every
$D\ge2$,
\begin{equation}\label{eq:main}
 c_{\Sigma,q}(1+\log D)^{|1/p-1/q|}
 \le \calS_p(E_D^{p,q};\Sigma)
 \le C_{\Sigma,q}(1+\log D)^{|1/p-1/q|}.
\end{equation}
In particular, for $X=L^{p,q}((0,1)^3)$,
\begin{equation}\label{eq:classification}
 \calS_p(X;\Sigma)<\infty\quad\Longleftrightarrow\quad q=p.
\end{equation}
For $p<q<\infty$, the space $X$ still has Gaussian type $p$.
\end{theorem}

The logarithm in \eqref{eq:main} is important. It gives not only a yes-or-no
answer for an infinite coefficient space but the exact rate at which a
finite collection of channels becomes unstable. At $q=2$ this rate is
\[
 \calS_{4/3}(\ell_D^{4/3,2};\Sigma)
 \asymp(1+\log D)^{1/4}.
\]
This is consistent with the positive \emph{scalar physical-domain}
$L^{4/3,2}\to L^2(\Sigma)$ restriction estimate: the two Lorentz spaces
play different roles. Theorem~\ref{thm:main} proves that attaining Gaussian
type $4/3$ does not by itself protect the critical Fourier trace.

\subsection{The fractional Gaussian comparison}

The same question can be asked for a sequence of correlated Gaussian
increments. Let $H=3/4$, and define the covariance of fractional Gaussian
noise by
\[
 r_H(m)=\frac12\Bigl(|m+1|^{2H}+|m-1|^{2H}-2|m|^{2H}\Bigr),
 \qquad F_N(j,k)=r_H(j-k).
\]
For deterministic coefficient vectors $x_j\in X$ and a centered
Gaussian vector $(Z_j)$ with covariance $F_N$, set
\begin{equation}\label{eq:fconst}
 \calF_p(X)=\sup_{N\ge1}\sup_{(x_j)\ne0}
 \frac{\bigl(\E\|\sum_{j=1}^N Z_jx_j\|_X^2\bigr)^{1/2}}
 {\bigl(\sum_{j=1}^N\|x_j\|_X^p\bigr)^{1/p}}.
\end{equation}
Unlike independent noise, these increments retain correlations over long
time intervals. The same Lorentz transition appears.

\begin{theorem}[Matching fractional Gaussian threshold]\label{thm:fractional}
For $p=4/3$ and $1\le q<\infty$,
\begin{equation}\label{eq:fractional-rate}
 \calF_p(E_D^{p,q})\asymp_q
 (1+\log D)^{|1/p-1/q|},\qquad D\ge2.
\end{equation}
Hence $\calF_p(L^{p,q}((0,1)^3))<\infty$ if and only if $q=p$.
\end{theorem}

The matching rates reveal a common difficulty in collecting many
correlated channels. They do \emph{not} prove that spherical restriction
and the fractional Gaussian model are equivalent for arbitrary Banach spaces.

\subsection{Results, a dispersive consequence and proof guide}

Theorem~\ref{thm:main} determines the precise cost of combining
many spectral channels at the critical Fourier exponent. Its lower bound
uses packets at successively larger parabolic scales: although the input
remains inexpensive in a Lorentz norm, each channel retains a nontrivial
Fourier trace. Theorem~\ref{thm:fractional} identifies the corresponding
threshold for correlated Gaussian increments, by a separate interval
construction. A further example shows that optimal Gaussian type and
cotype need not determine the restriction estimate.

There is also a direct dispersive PDE consequence. The linear capillary-wave
propagator has frequency law $\omega(\xi)=c|\xi|^{3/2}$, whose graph
is curved away from zero frequency. Theorem~\ref{thm:capillary} gives the
sharp dimension dependence of a spacetime $L^4$ estimate for arrays of
these waves, measured in Lorentz sequence norms. This application concerns
the linearized capillary-wave equation and requires no matrix-operator
model; the water-wave dispersion setting is discussed in
\cite{DengIonescuPausaderPusateri}.

The proof begins with finite-dimensional Lorentz inequalities and the
necessary Gaussian-type condition. It then treats the two sides of the
critical Lorentz index separately: multiscale Fourier packets for $q>4/3$
and Gaussian type for $q<4/3$. The fractional covariance argument follows
independently. Finally, we prove the capillary-wave consequence and record
an interpolation stability principle for the supercritical Lorentz range.

\section{Elementary tools: Lorentz norms, Gaussian type and scalarization}

\begin{lemma}[Finite Lorentz distortion]\label{lem:distortion}
For fixed $1\le q<\infty$ and $D\ge2$, the identity between $E_D^{p,q}$ and $\ell_D^p$ has Banach--Mazur distortion at most
\[
 C_{p,q}(1+\log D)^{|1/p-1/q|}.
\]
Moreover, $\norm{\one_D}_{E_D^{p,q}}\asymp_{p,q}D^{1/p}$, uniformly in $D$.
\end{lemma}
\begin{proof}
Let $z^*$ be nonincreasing, and partition the ranks into the dyadic
intervals $I_j=[2^j,2^{j+1})\cap\{1,\ldots,D\}$. There are
$J+1\asymp1+\log D$ intervals. Because $z^*$ is monotone, the
contribution of a full block to $\sum_r (z^*(r))^p$ is comparable to
$2^j(z^*(2^j))^p$. Similarly, the block contribution to the Lorentz
quasi-norm is comparable to $2^{jq/p}(z^*(2^j))^q$. Using a one-sided
adjacent-block shift at the final short block makes these comparisons
uniform. With
\[
 a_j=2^{j/p}z^*(2^j),
\]
this gives
\[
 \|z\|_{\ell^p}\asymp_p\Big(\sum_{j=0}^Ja_j^p\Big)^{1/p},
 \qquad
 \|z\|_{p,q}^*\asymp_{p,q}\Big(\sum_{j=0}^Ja_j^q\Big)^{1/q}.
\]
There is now no Fourier analysis in the calculation: it is the ordinary
comparison between $\ell^p$ and $\ell^q$ on a sequence of length $J+1$.
Its optimal loss is at most $(J+1)^{|1/p-1/q|}$. Replacing the Lorentz
quasi-norm by a uniformly equivalent lattice norm changes only a constant
that depends on $p,q$, not on $D$. Finally,
\[
 \|\one_D\|_{p,q}^{*q}
 =\sum_{r=1}^D r^{q/p-1}\asymp_{p,q} D^{q/p},
\]
so $\|\one_D\|_{E_D^{p,q}}\asymp D^{1/p}$.
\end{proof}

\begin{lemma}[Scalarization and the necessary type bound]\label{lem:type}
For $p=4/3$,
\[
 \calS_p(\ell_D^p;\Sigma)\le C_\Sigma,
 \qquad T_p(X)\le C_\Sigma\calS_p(X;\Sigma)
\]
for every Banach space $X$, whenever the right side is finite. In particular,
\begin{equation}\label{eq:dist-upper}
 \calS_p(E_D^{p,q};\Sigma)
 \lesssim_{\Sigma,q}(1+\log D)^{|1/p-1/q|}.
\end{equation}
\end{lemma}
\begin{proof}
For $X=\ell_D^p$, let $f=(f_1,\ldots,f_D)$. The Gaussian
$L^2$ norm in the definition of $\gamma$, followed by Gaussian moment
comparison and Fubini, gives
\[
 \|R_\Sigma f\|_{\gamma(L^2(\Sigma),\ell_D^p)}
 \asymp_p
 \left(\sum_{a=1}^D\|R_\Sigma f_a\|_2^p\right)^{1/p}.
\]
Apply scalar Stein--Tomas to every $f_a$ and note that
$\|f\|_{L^p(\R^3;\ell_D^p)}^p=\sum_a\|f_a\|_p^p$.
The resulting constant does not depend on $D$.

For necessity, fix a compactly supported scalar bump $\varphi$ whose
Fourier transform does not vanish on a positive-measure part of $\Sigma$.
Given $x_1,\ldots,x_m\in X$, choose translations $a_j$ successively so
far apart that the physical supports of the translated bumps are disjoint
and the Gram matrix of
\[
 h_j(\omega)=\widehat\varphi(\omega)e^{-ia_j\cdot\omega}
\]
has a lower spectral bound independent of $m$. To justify the latter,
observe that each off-diagonal Gram entry is the Fourier transform of the
smooth finite measure $|\widehat\varphi|^2d\sigma$ evaluated at
$a_j-a_k$. Those entries tend to zero as the distance grows. By making
them small enough that the sum of absolute off-diagonal entries in every
row is at most half the common diagonal entry, the Gram matrix is
bounded below by a fixed multiple of the identity.

Set $f(x)=\sum_{j=1}^m\varphi(x-a_j)x_j$. Disjoint physical supports give
\[
 \|f\|_{L^p(\R^3;X)}
 =\|\varphi\|_p\left(\sum_j\|x_j\|_X^p\right)^{1/p}.
\]
The lower Gram bound means that the Gaussian covariance of the restricted
family $(h_j)$ dominates a fixed multiple of the identity. Gaussian
covariance domination therefore yields
\[
 \|R_\Sigma f\|_{\gamma(L^2(\Sigma),X)}
 \gtrsim\left(\E\left\|\sum_j\gamma_jx_j\right\|_X^2\right)^{1/2}.
\]
The assumed restriction estimate gives precisely the type-$p$ inequality.
Finally, transfer the positive $\ell_D^p$ estimate to $E_D^{p,q}$ through
Lemma~\ref{lem:distortion}, paying the product of the norms of the identity
and its inverse.
\end{proof}

\begin{lemma}[Gaussian square-function lower bound]\label{lem:gamma-lower}
Let $E$ be a finite-dimensional Banach lattice with a lattice norm, and let $G:\Sigma\to E$ have finitely many scalar coordinates $G_a$. Then
\begin{equation}\label{eq:gamma-lower}
 \norm{G}_{\ga(L^2(\Sigma),E)}
 \ge c\left\|\left(\norm{G_a}_{L^2(\Sigma)}\right)_a\right\|_E
\end{equation}
with a numerical $c>0$. An equivalent lattice norm changes $c$ only by the equivalence constant.
\end{lemma}
\begin{proof}
The Gaussian radonification of $G$ has each coordinate a centered Gaussian of variance $\norm{G_a}_2^2$. Thus, using Jensen and the lattice order,
\[
 \norm{G}_\ga\ge\E\norm{Z}_E
 \ge\norm{(\E|Z_a|)_a}_E
 =\sqrt{2/\pi}\norm{(\norm{G_a}_2)_a}_E.
\]
For complex Gaussians the analogous fixed absolute-moment constant is used.
\end{proof}

\section{Multiscale wave packets: the case \texorpdfstring{$q>p$}{q>p}}

This section proves the main lower bound for $p<q<\infty$ without
using scalar Lorentz sharpness. The proof has three ingredients: (i) one
packet at any parabolic scale has bounded $L^p$ input and nontrivial
restriction output; (ii) the same physical point meets only a controlled
number of translated packets at each scale; and (iii) the restriction
outputs retain a separate coordinate for \emph{every} translated packet.
The last two facts are what create the sharp logarithmic gap.

Fix orthonormal physical coordinates $x=(x_\perp,x_3)$ adapted to a patch of $\Sigma$ near $\xi_0$; locally $\Sigma$ is a graph $\xi_3=\xi_{0,3}+O(|\xi_\perp-\xi_{0,\perp}|^2)$. Choose $\Psi\in C_c^\infty(\R^3)$ such that $\widehat\Psi(0)\ne0$, and set for $R\ge1$
\begin{equation}\label{eq:packet}
 g_R(x)=R^{-3}e^{ix\cdot\xi_0}
 \Psi\!\left(\frac{x_1}{R},\frac{x_2}{R},\frac{x_3}{R^2}\right).
\end{equation}
Its $L^p$ norm is uniformly bounded, since the physical support volume is $O(R^4)$ and $R^{-3p}R^4=1$ at $p=4/3$.

\begin{lemma}[Uniform scalar packet output]\label{lem:packet}
There is $c_{\Sigma,\Psi}>0$ with
\[
 \norm{R_\Sigma g_R}_{L^2(\Sigma)}\ge c_{\Sigma,\Psi},
 \qquad R\ge1.
\]
\end{lemma}
\begin{proof}
The change of variables $x_1=Ry_1$, $x_2=Ry_2$,
$x_3=R^2y_3$ has Jacobian $R^4$. Hence
\[
 \widehat g_R(\xi)=R\,\widehat\Psi\bigl(
 R(\xi_\perp-\xi_{0,\perp}),
 R^2(\xi_3-\xi_{0,3})\bigr).
\]
Write the local surface graph as
$\xi_3=\xi_{0,3}+\phi(\xi_\perp-\xi_{0,\perp})$, where
$\phi(0)=0$, $\nabla\phi(0)=0$ and
$|\phi(y)|\le C|y|^2$ for small $y$. Since $\widehat\Psi(0)\ne0$,
its absolute value is bounded below in a sufficiently small neighborhood
of the origin. On the cap $|y|\le c_0/R$, both
$Ry$ and $R^2\phi(y)$ lie in this neighborhood. Therefore
$|\widehat g_R|\gtrsim R$ throughout that cap. Its surface measure is
comparable to $R^{-2}$, so
\[
 \|R_\Sigma g_R\|_2^2
 \ge\int_{|y|\le c_0/R}|\widehat g_R|^2\,d\sigma
 \gtrsim R^2R^{-2}\gtrsim1.
\]
The finitely bounded range of smaller $R$ is handled by reducing the
constant after choosing the patch.
\end{proof}

Take $M=2^h$, $h\ge1$, and put
\[
 N=M^4=2^{4h},\qquad J=4h,\qquad R_j=2^{j/4}
 \quad(0\le j<J).
\]
At level $j$ the packet occupies a physical box of volume $R_j^4=2^j$
and has height $R_j^{-3}=2^{-j/p}$. The choice $N=2^J$ allows us to fit
all $J$ levels in one common large box. This balance is the source of the
logarithm: a level contributes about $2^j$ nonzero coordinates, but each is
only about $2^{-j/p}$ in magnitude.
Choose the lattice of centers
\[
 \mathcal A_M=\{1,\dots,M\}^2\times\{1,\dots,M^2\},
 \qquad |\mathcal A_M|=N,
\]
expressed in the adapted physical coordinates. The coefficient index set is
$\mathcal I=\{0,\dots,J-1\}\times\mathcal A_M$, of cardinality $D=NJ$.
Define $\mathcal G:\R^3\to E_D^{p,q}$ by
\begin{equation}\label{eq:array}
 \mathcal G(x)_{(j,a)}=g_{R_j}(x-a).
\end{equation}

\begin{lemma}[Pointwise weak-$\ell^p$ envelope]\label{lem:envelope}
The array \eqref{eq:array} satisfies
\[
 \sup_{x\in\R^3}\norm{\mathcal G(x)}_{\ell_D^{p,\infty}}\le C_\Psi,
 \quad |\operatorname{supp}\mathcal G|\le C_\Psi N,
\]
and hence, for each fixed finite $q$,
\begin{equation}\label{eq:input}
 \norm{\mathcal G}_{L^p(\R^3;E_D^{p,q})}
 \le C_{\Psi,p,q}N^{1/p}J^{1/q}.
\end{equation}
\end{lemma}
\begin{proof}
At scale $j$, every nonzero coordinate has absolute value at most $C R_j^{-3}=C2^{-j/p}$. Since $\Psi$ is compactly supported, at a fixed $x$ there are at most
$C(R_j+1)^2(R_j^2+1)\le C R_j^4=C2^j$
centers $a$ with $g_{R_j}(x-a)\ne0$. Consequently, for $t>0$,
\[
 \#\{(j,a):|\mathcal G(x)_{(j,a)}|>t\}
 \le C\sum_{2^{-j/p}\gtrsim t}2^j
 \le C't^{-p}.
\]
The key point is that the count grows geometrically: the entire sum
through the last active level is comparable to its final term. Thus the
counting bound is $Ct^{-p}$, independent of the number of scales.
This is the uniform weak-$\ell^p$ estimate. Also $R_j\le M$ and
$R_j^2\le M^2$, so all the translated supports lie in a fixed enlargement
of the box $M\times M\times M^2$, whose volume is $O(N)$.

To pass from weak-$\ell^p$ to the stronger $\ell^{p,q}$ norm, order the
nonzero coordinates in decreasing magnitude. The weak bound says that the
$r$-th largest coordinate is at most $Cr^{-1/p}$. Substitution into
\eqref{eq:lorentz} gives
\[
 \norm{z}_{p,q}^q\lesssim\sum_{r=1}^D r^{q/p-1}(Cr^{-1/p})^q
 \lesssim 1+\log D.
\]
Since $D=NJ$ and $N=2^J$, $1+\log D\asymp J$. Integration over the union of the supports proves \eqref{eq:input}.
\end{proof}

\begin{proof}[Proof of the lower bound in Theorem~\ref{thm:main} for $q>p$]
Every scalar coordinate of $R_\Sigma\mathcal G$ is a modulation of $\widehat g_{R_j}|_\Sigma$ by $e^{-ia\cdot\xi}$, hence has $L^2(\Sigma)$ norm at least $c_{\Sigma,\Psi}$. Lemma~\ref{lem:gamma-lower} therefore gives
\[
 \norm{R_\Sigma\mathcal G}_{\ga(L^2(\Sigma),E_D^{p,q})}
 \gtrsim_{p,q,\Sigma}\norm{\one_D}_{E_D^{p,q}}
 \asymp D^{1/p}=(NJ)^{1/p}.
\]
Divide by \eqref{eq:input} to obtain
\[
 \calS_p(E_D^{p,q};\Sigma)\gtrsim J^{1/p-1/q}.
\]
Here $D=2^JJ$, so $J\asymp1+\log D$. The restricted sequence of
values $J=4h$ still gives coefficient dimensions with bounded successive
ratios. Given an intermediate $D$, choose the largest constructed dimension
not exceeding $D$ and append zero coordinates; the norm of this inclusion
is controlled independently of dimension. This extends the lower bound to
all large $D$, and the remaining finite set of dimensions is covered by
rank-one scalar tests.
\end{proof}

\begin{remark}[Where the sharp loss occurs]
At every physical point the multiscale array has a controlled
weak-$\ell^p$ profile, so the source pays only $J^{1/q}$ for its $J$
levels. By contrast, after restriction each one of the $D=NJ$ coordinates
has an $L^2(\Sigma)$ norm bounded below; Gaussian mixing then pays
$J^{1/p}$. The quotient $J^{1/p-1/q}$ is therefore forced by the difference
between physical aggregation and Gaussian output aggregation. It is not a
loss in the scalar Stein--Tomas estimate for a single packet.
\end{remark}

\subsection{The same obstruction in higher dimensions}

The physical-counting mechanism is not peculiar to the three-dimensional sphere. Its decisive feature is that the Knapp packet support volume and the number of overlapping packet centers have the same parabolic exponent.

\begin{corollary}[Curved hypersurfaces in every ambient dimension]\label{cor:higher}
Let $d\ge2$ and $p_d=2(d+1)/(d+3)$. Suppose $\Sigma\subset\R^d$ is a compact smooth hypersurface with nonvanishing Gaussian curvature and a smooth positive density. For every fixed $1\le q<\infty$ and $D\ge2$,
\[
 \norm{R_\Sigma:L^{p_d}(\R^d;\ell_D^{p_d,q})
 \to\gamma(L^2(\Sigma),\ell_D^{p_d,q})}
 \asymp_{d,q,\Sigma}(1+\log D)^{|1/p_d-1/q|}.
\]
Consequently $L^{p_d,q}$ as a coefficient space satisfies the critical $\gamma$-valued spherical restriction bound if and only if $q=p_d$.
\end{corollary}
\begin{proof}
The scalar Stein--Tomas theorem and Gaussian Fubini give the upper bound on $\ell_D^{p_d}$, and finite Lorentz distortion gives the stated upper bound. The necessary Gaussian type-$p_d$ argument and disjoint-block construction from Lemmas~\ref{lem:type} and \ref{lem:below} work in every dimension, proving the lower bound when $q<p_d$.

For $q>p_d$, near a point of $\Sigma$ with tangent-normal coordinates $x=(x'\!,x_d)$, use
\[
 g_R(x)=R^{-(d+1)/p_d}e^{ix\cdot\xi_0}
 \Psi(x'/R,x_d/R^2).
\]
Its physical support has volume $O(R^{d+1})$, and its Fourier transform has magnitude comparable to
$R^{(d+1)(1-1/p_d)}=R^{(d-1)/2}$ on a cap of surface measure $\gtrsim R^{-(d-1)}$. Thus $\norm{R_\Sigma g_R}_2\gtrsim1$.
Choose $M=2^h$, $N=M^{d+1}$, $J=(d+1)h$, and $R_j=2^{j/(d+1)}$. At each scale translate $g_{R_j}$ over a common grid
$\{1,\ldots,M\}^{d-1}\times\{1,\ldots,M^2\}$, containing $N$ centers. Every point belongs to at most $C R_j^{d+1}=C2^j$ of the supports at scale $j$, and the amplitude of each such coordinate is at most $C2^{-j/p_d}$. The coefficient array therefore has pointwise weak-$\ell^{p_d}$ norm $O(1)$, Lorentz $\ell^{p_d,q}$ norm $O(J^{1/q})$, and spatial support volume $O(N)$. Its input is at most $CN^{1/p_d}J^{1/q}$, while Lemma~\ref{lem:gamma-lower} gives output at least $c(NJ)^{1/p_d}$. This proves the required lower rate, and padding finishes the proof.
\end{proof}

\section{The Gaussian-type obstruction when \texorpdfstring{$q<p$}{q<p}}

For $q<p$ the obstruction is already present at the Gaussian-type level; curvature-specific packets are unnecessary.

\begin{lemma}[Critical type failure below the Lorentz diagonal]\label{lem:below}
Let $1\le q<p$. For $D\ge2$,
\[
 T_p(E_D^{p,q})\gtrsim_{p,q}(1+\log D)^{1/q-1/p}.
\]
\end{lemma}
\begin{proof}
Choose $J\ge2$ and pairwise disjoint coordinate blocks $B_j$ with
$|B_j|=2^j$, $0\le j<J$. Put $y_j=2^{-j/p}\one_{B_j}$.
The decreasing rearrangement of $y_j$ is constant on its nonzero block,
and a direct substitution in \eqref{eq:lorentz} shows
$\|y_j\|_{p,q}\asymp1$ uniformly in $j$.

The sum $y=\sum_jy_j$ has decreasing values approximately $r^{-1/p}$
for ranks $1\le r\lesssim2^J$. Each dyadic block consequently contributes
a fixed positive amount to the $q$-th power of the Lorentz norm, giving
\[
 \|y\|_{p,q}\asymp J^{1/q}.
\]
Because the supports are disjoint, the absolute value of the Gaussian sum
is the vector $\sum_j|\gamma_j|y_j$. Jensen's inequality for the lattice
norm and the positive cone yields
\[
 \begin{aligned}
 \left(\E\left\|\sum_j\gamma_jy_j\right\|_{p,q}^2\right)^{1/2}
 &\ge \E\left\|\sum_j\gamma_jy_j\right\|_{p,q}\\
 &\ge \left\|\sum_j(\E|\gamma_j|)y_j\right\|_{p,q}
 \gtrsim J^{1/q}.
 \end{aligned}
\]
The type-$p$ denominator is
$\bigl(\sum_j\|y_j\|_{p,q}^p\bigr)^{1/p}\asymp J^{1/p}$.
Finally, $D=\sum_{j<J}2^j\asymp2^J$; appending zero coordinates covers
other sufficiently large dimensions.
\end{proof}

Lemma~\ref{lem:type} and Lemma~\ref{lem:below} give the lower spherical bound for $q<p$. The upper bound for every $q$ follows from Lemmas~\ref{lem:distortion}--\ref{lem:type}. This proves \eqref{eq:main}.

The infinite-dimensional statement \eqref{eq:classification} follows because $E_D^{p,q}$ embeds (with the standard counting-measure rescaling) as a sublattice of $L^{p,q}((0,1)^3)$ with constants independent of $D$. For $q=p$, $L^{p,p}=L^p$, and Gaussian Fubini plus scalar Stein--Tomas give the positive endpoint. For $q>p$, the lower bound just proved diverges; for $q<p$, the type obstruction diverges. Creekmore's theorem \cite{Creekmore} supplies attained Gaussian type $p$ for the infinite Lorentz spaces $L^{p,q}$ with $q>p$ and finite $q$. This completes Theorem~\ref{thm:main}.

\section{The fractional comparison model}

We record the matching Lorentz result for the covariance \eqref{eq:fconst}, keeping its proof independent of the spherical construction.

\begin{lemma}[Scalar fractional amplification]\label{lem:scalarF}
For $p=4/3$, $\calF_p(\ell_D^p)\le C$ independently of $D$. Moreover, $T_p(X)\lesssim\calF_p(X)$ for all Banach spaces $X$.
\end{lemma}
\begin{proof}
A Taylor expansion of $|m\pm1|^{2H}$ at large $|m|$ gives
\[
 r_H(m)=H(2H-1)|m|^{2H-2}+O(|m|^{2H-4}).
\]
At $H=3/4$ this is
$r_H(m)=\tfrac38|m|^{-1/2}+O(|m|^{-5/2})$, and in particular
$|r_H(m)|\le C(1+|m|)^{-1/2}$. For a finitely supported scalar
sequence $a=(a_j)$, split the covariance form into the diagonal and
off-diagonal parts. The diagonal is bounded by $\|a\|_{\ell^p}^2$
since $p<2$. Discrete Hardy--Littlewood--Sobolev controls the off-diagonal
part:
\[
 \sum_{j\ne k}\frac{|a_j|\,|a_k|}{|j-k|^{1/2}}
 \le C\|a\|_{\ell^{4/3}}^2.
\]
Therefore the scalar Gaussian covariance is bounded by
$C\|a\|_{\ell^p}^2$, uniformly in the length of the sequence.
Gaussian moment comparison and coordinatewise Fubini now give
$\calF_p(\ell_D^p)\le C$ independently of $D$.

To see necessity of type, take any fixed number $m$ of coefficient vectors,
place them at widely separated indices in one sufficiently long window,
and put zero vectors at all other indices. Since $r_H(n)\to0$, the resulting
$m\times m$ covariance matrix can be made arbitrarily close to the
identity in operator norm. Gaussian covariance domination then recovers
an independent Gaussian sum, proving $T_p(X)\lesssim\calF_p(X)$.
\end{proof}

\begin{proof}[Proof of Theorem~\ref{thm:fractional}]
We first separate the easy bounds from the multiscale lower bound. The
upper estimate for every $q$ follows by applying Lemma~\ref{lem:distortion}
to the scalar fractional estimate of Lemma~\ref{lem:scalarF}.
For $q<p$, combine the necessary type inequality in
Lemma~\ref{lem:scalarF} with the disjoint-block example of
Lemma~\ref{lem:below}. The case $q=p$ is uniformly bounded by the scalar
estimate and Gaussian Fubini.

It remains to prove the lower bound for $q>p$. Fix $J\ge2$ and let
$N=2^J$. At level $j=0,\ldots,J-1$, use an interval length $L_j=2^j$.
Give one separate coefficient coordinate to every interval
$[z,z+L_j-1]\subset\{1,\ldots,N\}$. Denote its label by $(j,z)$ and let
\[
 \mathcal I=\{(j,z):0\le j<J,\;1\le z\le N-L_j+1\}.
\]
The total number of coefficient coordinates is
\[
 D=|\mathcal I|=\sum_{j=0}^{J-1}(N-L_j+1)\asymp NJ.
\]
For each time index $1\le n\le N$, define the coefficient vector
$x_n\in E_D^{p,q}$ by
\begin{equation}\label{eq:fractional-test-input}
 x_n(j,z)=L_j^{-1/p}\one_{\{z\le n<z+L_j\}}.
\end{equation}

\emph{Input counting.} At level $j$, at most $L_j=2^j$ intervals
contain a given $n$, and every corresponding coordinate of $x_n$ has
size $L_j^{-1/p}=2^{-j/p}$. Exactly as in the spherical packet-counting
lemma, the number of coordinates greater than $t$ is $O(t^{-p})$.
Therefore $\|x_n\|_{\ell^{p,\infty}}\le C$ and
$\|x_n\|_{E_D^{p,q}}\lesssim J^{1/q}$. The total input budget is
\begin{equation}\label{eq:fractional-input-budget}
 \left(\sum_{n=1}^N\|x_n\|_{E_D^{p,q}}^p\right)^{1/p}
 \lesssim N^{1/p}J^{1/q}.
\end{equation}

\emph{Gaussian output.} Let $(Z_n)$ be fractional Gaussian noise with
covariance $r_H(j-k)$. By definition, there is a fractional Brownian
motion $B_H$ such that $Z_n=B_H(n)-B_H(n-1)$. Hence, for every interval
of $L$ consecutive increments,
\[
 \sum_{n=z}^{z+L-1} Z_n=B_H(z+L-1)-B_H(z-1),
 \qquad \mathrm{Var}\!\left(\sum_{n=z}^{z+L-1}Z_n\right)=L^{2H}.
\]
Since $H=1/p$, multiplying by the coefficient $L^{-1/p}$ makes this
variance exactly one. Thus \emph{each} coordinate of the vector
$\sum_nZ_nx_n\in E_D^{p,q}$ is a centered Gaussian with unit variance,
although the different coordinates are correlated. The Gaussian lattice
lower bound of Lemma~\ref{lem:gamma-lower} applies and gives
\begin{equation}\label{eq:fractional-output-budget}
 \left(\E\left\|\sum_nZ_nx_n\right\|_{E_D^{p,q}}^2\right)^{1/2}
 \gtrsim \|\one_D\|_{E_D^{p,q}}
 \asymp D^{1/p}\asymp(NJ)^{1/p}.
\end{equation}

Divide \eqref{eq:fractional-output-budget} by
\eqref{eq:fractional-input-budget}. The quotient is at least
$J^{1/p-1/q}$. Since $D\asymp 2^JJ$, we have
$J\asymp1+\log D$. Padding by zero coordinates extends the result to
all dimensions. This completes the fractional lower bound and the theorem.
\end{proof}

\section{A sharp Lorentz threshold for linear capillary waves}\label{sec:capillary}

The endpoint restriction theorem also gives a precise spacetime estimate
for a dispersive partial differential equation. The example is the
linearized, infinite-depth capillary-wave model on a two-dimensional water
surface. After diagonalizing the linearized system and neglecting gravity,
its two time-frequency branches have dispersion relation
\[
 \omega(\xi)=c|\xi|^{3/2},\qquad c>0,\quad \xi\in\R^2.
\]
The corresponding second-order equation is
\begin{equation}\label{eq:capillary-pde}
 \partial_t^2 U+c^2(-\Delta_x)^{3/2}U=0,
 \qquad(t,x)\in\R\times\R^2.
\end{equation}
Here $(-\Delta_x)^{3/2}$ is the Fourier multiplier with symbol $|\xi|^3$.
The fractional power is characteristic of capillary dispersion; see
\cite{DengIonescuPausaderPusateri} for the broader water-wave context.
We study one complex frequency branch, which is sufficient to exhibit
the sharp coefficient-space threshold.

Let $\mathcal A=\{\xi\in\R^2:1<|\xi|<2\}$ and fix $1<s<\infty$.
For $D$ initial profiles $U_0=(u_{0,1},\ldots,u_{0,D})$, set
$Y_D^{s}=\ell_D^{4,s}$, with a uniformly equivalent Banach lattice norm.
Write
\[
 U(t)=e^{-itc(-\Delta_x)^{3/4}}U_0,
 \qquad \operatorname{supp}\widehat{u}_{0,j}\Subset\mathcal A.
\]
Thus $U$ solves \eqref{eq:capillary-pde}, with
$\partial_tU(0)=-ic(-\Delta_x)^{3/4}U_0$.

\begin{theorem}[Sharp capillary-wave spacetime estimate]\label{thm:capillary}
For each $c>0$ and $1<s<\infty$, there exist constants
$0<c_{s,c}\le C_{s,c}<\infty$, independent of $D\ge2$, such that
\begin{equation}\label{eq:capillary-bound}
 \|U\|_{L^4(\R_t\times\R_x^2;\ell_D^{4,s})}
 \le C_{s,c}(1+\log D)^{|1/4-1/s|}
 \left\|\left(\|\widehat{u}_{0,j}\|_{L^2(\mathcal A)}\right)_{j=1}^D
 \right\|_{\ell_D^{4,s}}.
\end{equation}
The power of $1+\log D$ is sharp: the best constant in
\eqref{eq:capillary-bound}, over all smooth data whose Fourier
transforms are supported in $\mathcal A$, is bounded below by
$c_{s,c}(1+\log D)^{|1/4-1/s|}$.
In particular, the estimate is dimension-independent at $s=4$,
whereas at $s=2$ its optimal loss is $(1+\log D)^{1/4}$.
\end{theorem}

The index $s$ is a \emph{coefficient} Lorentz exponent: the time-space
norm remains the ordinary $L^4$ norm. Thus the theorem is not a
Lorentz refinement in time or in physical space. It determines how the
same dispersive waves behave when many scalar components are combined
with different sequence-space geometries.

We record the two functional-analytic ingredients used to pass from
restriction to the wave equation. They also explain why the logarithmic
lower bound transfers, rather than giving only an upper estimate.

\begin{lemma}[Gaussian norm and trace duality]\label{lem:capillary-duality}
Fix $1<s<\infty$ and write $s'=s/(s-1)$,
$X_D=\ell_D^{4/3,s'}$ and $Y_D=\ell_D^{4,s}$.
Let $H$ be a Hilbert space. Given $G_1,\ldots,G_D\in H$, regard
$G=(G_j)_{j=1}^D$ as the finite-rank operator
$T_G:H\to Y_D$, $T_Gh=(\langle h,G_j\rangle_H)_{j=1}^D$.
Then
\begin{equation}\label{eq:capillary-gamma-sq}
 \|G\|_{\gamma(H,Y_D)}
 \asymp_s
 \left\|\left(\|G_j\|_H\right)_{j=1}^D\right\|_{Y_D}.
\end{equation}
Moreover, for a compact curved surface $\Sigma$ with smooth positive
density and $H=L^2(\Sigma)$,
\begin{equation}\label{eq:capillary-trace-dual}
 \|E_\Sigma:\gamma(H,Y_D)\to L^4(\R^3;Y_D)\|
 \asymp_s
 \mathcal S_{4/3}(X_D;\Sigma),
\end{equation}
with constants independent of $D$.
\end{lemma}
\begin{proof}
Choose an orthonormal expansion
$G_j=\sum_{\nu=1}^m a_{\nu,j}h_\nu$ in $H$.
The associated $\gamma$-radonifying norm is the second Gaussian
moment of the $Y_D$-valued vector
$Z=(\sum_\nu\gamma_\nu a_{\nu,j})_{j=1}^D$.
The Gaussian square-function theorem for finite-cotype Banach lattices
states that this moment is comparable to the lattice norm of its
coordinate standard deviations; see \cite{vanNeerven}.
The standard deviations are
$(\sum_\nu|a_{\nu,j}|^2)^{1/2}=\|G_j\|_H$.
The Lorentz sequence spaces $\ell_D^{4,s}$ have finite-cotype constants
independent of $D$ for fixed $s$, so this proves
\eqref{eq:capillary-gamma-sq}.

For \eqref{eq:capillary-trace-dual}, Lorentz duality gives
$X_D^*\simeq Y_D$, with equivalence constants independent of $D$;
see \cite{BennettSharpley,Creekmore}.
The spaces $X_D$ are uniformly $K$-convex. Indeed, for fixed $s'$
they arise by real interpolation between $\ell_D^{r_0}$ and
$\ell_D^{r_1}$ for suitable $1<r_0<4/3<r_1<\infty$;
these interpolation families have nontrivial type, with constants
independent of $D$. The $\gamma$-trace duality theorem therefore identifies
$(\gamma(H,X_D))^*$ with $\gamma(H,X_D^*)$, up to a constant
independent of $D$; see \cite{vanNeerven}.
On finite-rank functions, the adjoint of Fourier restriction is the
extension operator:
\[
 \int_\Sigma\langle R_\Sigma f(\omega),G(\omega)\rangle\,d\sigma(\omega)
 =\int_{\R^3}\langle f(x),E_\Sigma G(x)\rangle\,dx,
\]
with the usual harmless Fourier-conjugation convention.
Bochner $L^{4/3}$--$L^4$ duality and the preceding $\gamma$-duality
identify their operator norms within the claimed uniform constants.
\end{proof}

\begin{proof}[Proof of Theorem~\ref{thm:capillary}]
The solution is the oscillatory integral
\begin{equation}\label{eq:capillary-ext}
 U(t,x)=(2\pi)^{-2}\int_{\mathcal A}
 e^{i(x\cdot\xi-tc|\xi|^{3/2})}
 \widehat U_0(\xi)\,d\xi.
\end{equation}
It is a Fourier extension integral for the graph
\[
 \Sigma_c=\{(\xi,c|\xi|^{3/2}):\xi\in\mathcal A\}
 \subset\R^3.
\]
For a radial graph $\tau=\omega(|\xi|)$, the Hessian in the
frequency plane has eigenvalues $\omega''(r)$ in the radial direction
and $\omega'(r)/r$ in the tangential direction.
For $\omega(r)=cr^{3/2}$ these are
\[
 \omega''(r)=\frac{3c}{4}r^{-1/2},
 \qquad \frac{\omega'(r)}r=\frac{3c}{2}r^{-1/2},
\]
so
\[
 \det D^2\omega(\xi)=\frac{9c^2}{8|\xi|}>0
 \quad\text{on }\mathcal A.
\]
Thus $\Sigma_c$ has nonvanishing Gaussian curvature on the
chosen annulus. The graph measure $d\xi$ differs from surface
measure by a bounded, smooth positive factor, which changes only
constants depending on $c$ and the annulus.

Apply Theorem~\ref{thm:main} with Lorentz secondary index $q=s'$:
\[
 \mathcal S_{4/3}(X_D;\Sigma_c)
 \asymp_{s,c}(1+\log D)^{|3/4-1/s'|}
 =(1+\log D)^{|1/4-1/s|}.
\]
Lemma~\ref{lem:capillary-duality} transfers this sharp estimate to
$E_{\Sigma_c}: \gamma(L^2(\mathcal A),Y_D)\to
L^4(\R^{1+2};Y_D)$.
Finally, formula \eqref{eq:capillary-gamma-sq} turns the
$\gamma$-norm of the frequency data into the coefficient norm on
the right of \eqref{eq:capillary-bound}. This proves its upper bound.
The lower bound is equally important: if a smaller power of
$1+\log D$ held uniformly for every smooth initial array,
\eqref{eq:capillary-trace-dual} would imply a smaller restriction
constant for $X_D$, contradicting the lower half of
Theorem~\ref{thm:main}. Smooth compactly supported frequency
arrays are dense in the finite-rank spaces used in the duality argument,
so the sharpness is realized up to a fixed multiplicative constant
by such arrays.\end{proof}

\begin{remark}[Why the coefficient threshold is relevant]
The result distinguishes two notions of stability for a family of
linear waves. Each scalar capillary mode obeys the familiar $L^4$
spacetime estimate, but assembling $D$ modes in $\ell^{4,s}$ can
cost an unavoidable logarithmic power. At $s=4$ the channels combine
uniformly; at $s=2$ a fourth-root logarithm is necessary.
The statement concerns the linearized model, not the nonlinear
free-boundary evolution.
\end{remark}

\section{Real interpolation beyond the critical Lorentz line}\label{sec:interpolation}

The sharp logarithmic estimate of Theorem~\ref{thm:main} concerns the \emph{critical primary} Lorentz exponent $p=4/3$. A complementary positive region for larger primary exponents follows from a general stability principle. The argument is included to distinguish a genuine endpoint obstruction from ordinary interpolation of positive spaces. We use standard real-interpolation conventions as in \cite{Interpolation,BennettSharpley}.

\begin{proposition}[Real-interpolation stability of critical restriction]\label{prop:interpolation}
Let $(X_0,X_1)$ be a compatible regular Banach couple, let $0<\theta<1$, and put $X_{\theta,q}=(X_0,X_1)_{\theta,q}$ using the real $K$-method. If
$R_\Sigma:L^p(\R^3;X_i)\to\gamma(L^2(\Sigma),X_i)$ is bounded for $i=0,1$ and $p=4/3$, then for every finite $p\le q<\infty$,
\[
 \calS_p(X_{\theta,q};\Sigma)
 \le C_{p,q,\theta}\,\calS_p(X_0;\Sigma)^{1-\theta}
                 \calS_p(X_1;\Sigma)^\theta.
\]
\end{proposition}
\begin{proof}
Let $H=L^2(\Sigma)$. For the source, the $p$-convexified $K$-functional
\[
 K_p(t,x;X_0,X_1)=\inf_{x=x_0+x_1}
       (\|x_0\|_{X_0}^p+t^p\|x_1\|_{X_1}^p)^{1/p}
\]
is uniformly equivalent to the usual $K$-functional. On finite simple Bochner functions the $K_p$-functional for the couple $(L^p(X_0),L^p(X_1))$ is the $L^p$ norm of the pointwise $K_p$-functional. The standard density and measurability argument extends this comparison to the compatible Bochner couple. Minkowski for $q\ge p$ therefore gives
\[
 L^p(\R^3;(X_0,X_1)_{\theta,q})
 \hookrightarrow (L^p(\R^3;X_0),L^p(\R^3;X_1))_{\theta,q}.
\tag{A}
\]
For the target, start with a finite-rank $T:H\to X_0+X_1$ and an orthonormal basis of its Hilbert input. Let $Z_T=\sum_j\gamma_jTh_j$. For any decomposition $T=T_0+T_1$ and almost every Gaussian outcome,
\[
 K(t,Z_T;X_0,X_1)\le \|Z_{T_0}\|_{X_0}+t\|Z_{T_1}\|_{X_1}.
\]
Taking $L^q$ in the Gaussian variable, then infimizing over operator decompositions and using Gaussian Kahane--Khintchine, yields
\[
 \|K(t,Z_T;X_0,X_1)\|_{L^q(\Omega)}
 \le C_q K(t,T;\gamma(H,X_0),\gamma(H,X_1)).
\]
Integrating its $q$-th power against $t^{-\theta q}\,dt/t$, applying Fubini, and again using Kahane--Khintchine gives
\[
 (\gamma(H,X_0),\gamma(H,X_1))_{\theta,q}
 \hookrightarrow\gamma(H,(X_0,X_1)_{\theta,q}).
\tag{B}
\]
These estimates extend from finite-rank operators by completion. Interpolating the same restriction operator between its two endpoint bounds and composing (A)--(B) proves the assertion.
\end{proof}

\begin{corollary}[Positive Lorentz region above the critical primary exponent]\label{cor:supercritical}
For $r>4/3$ and $4/3\le q<\infty$, every classical Lorentz
function space $L^{r,q}$ satisfies the critical Stein--Tomas restriction
estimate
\[
 \mathcal S_{4/3}(L^{r,q};\Sigma)<\infty.
\]
\end{corollary}
\begin{proof}
Choose $4/3\le r_0<r<r_1<\infty$ and $0<\theta<1$ with
$1/r=(1-\theta)/r_0+\theta/r_1$. The real-interpolation identity
$(L^{r_0},L^{r_1})_{\theta,q}=L^{r,q}$ holds with equivalence
constants independent of finite-dimensional approximations.
For $r_i\ge4/3$, scalar Stein--Tomas and the mixed-norm
Minkowski inequality yield the critical estimate on $L^{r_i}$:
the $L^2$ norm in surface frequency can be passed inside the
coefficient $L^{r_i}$ norm, and the physical $L^{4/3}$ norm controls
the resulting mixed integral since $r_i\ge4/3$.
Apply Proposition~\ref{prop:interpolation} to conclude.
\end{proof}

\begin{corollary}[Matching type and cotype need not predict restriction]\label{cor:matching}
Fix $q>2$ and put
\[
 X_-=L^{4/3,q}(0,1),\qquad
 X_+=L^{4/3}\big((0,1);L^q(0,1)\big).
\]
Both spaces have optimal attained Gaussian type $4/3$ and optimal cotype $q$, but
\[
 \calS_{4/3}(X_-;S^2)=\infty,\qquad
 \calS_{4/3}(X_+;S^2)<\infty.
\]
\end{corollary}
\begin{proof}
The type and cotype exponents follow from the classical Lorentz calculations \cite{Creekmore} and the standard mixed-Lebesgue geometry. The negative assertion is Theorem~\ref{thm:main}. For the positive assertion, $L^q(0,1)$ has actual Gaussian type $2$ and therefore critical restriction; outer $L^{4/3}$-integration preserves the estimate by Gaussian Fubini and Minkowski.
\end{proof}

The index restriction $q<\infty$ in Proposition~\ref{prop:interpolation} is essential for the argument: Gaussian $L^q$-moment interchange cannot simply be replaced by an essential supremum at $q=\infty$. In the usual infinite-dimensional weak-Lorentz spaces $L^{r,\infty}$ for $r>1$, disjoint blocks of sizes $2^j$, scaled by $2^{-j/r}$, span uniformly isomorphic copies of $\ell^\infty_n$. Hence these weak-Lorentz coefficient spaces fail the necessary nontrivial Gaussian type condition; the same obstruction passes to the diagonal of $S^{r,\infty}$. These corollaries do not settle the general Banach-space classification at attained type $4/3$.

\section{Geometric interpretation and limitations}

Theorems~\ref{thm:main} and \ref{thm:fractional} identify a precise
boundary between independent Gaussian aggregation and the accumulation
of multiscale correlations. At the critical primary Lorentz exponent,
a finite-dimensional distortion of logarithmic size can become visible
in the Fourier trace even though every scalar component satisfies the
same Stein--Tomas estimate.

Gaussian type is related to Besov-to-$\gamma$ embeddings by
Kalton--van Neerven--Veraar--Weis \cite{KNVW}. Critical
Bessel-potential-to-$\gamma$ embeddings and positive lattice cases
are studied in \cite{Veraar}. These provide useful context, but type
and cotype alone do not characterize the curved restriction endpoint:
Corollary~\ref{cor:matching} supplies two explicit spaces with matching
optimal type and cotype and different restriction behavior.

The matching Lorentz growth laws do not imply either universal Banach
implication
\[
 \mathcal S_{4/3}(X)<\infty\ \Longrightarrow\
 \mathcal F_{4/3}(X)<\infty,\qquad
 \mathcal F_{4/3}(X)<\infty\ \Longrightarrow\
 \mathcal S_{4/3}(X)<\infty.
\]
Both would require estimates uniform across all finite-dimensional
coefficient geometries, beyond the Lorentz family. The
operator/Banach-space bipolar theorem of de la Salle \cite{delaSalle}
provides an abstract framework for such universal implications;
it does not itself supply a comparison for these two operators.

\section{The critical Gaussian-type threshold and open problems}

The results above separate Gaussian type from the curvature-sensitive endpoint.
At $p=4/3$, Gaussian type $p$ is necessary but is not sufficient for either
$\calS_p$ or $\calF_p$, while strict Gaussian type $r>p$ is known to imply
the critical fractional $\gamma$-embedding by the type--$\gamma$ theory of
\cite{Veraar}. This leaves a genuine comparison question at attained type $p$.
It is important not to conjecture that type $p$ alone suffices: Theorems
\ref{thm:main} and \ref{thm:fractional} disprove that statement.

\begin{conjecture}[Critical spherical-to-fractional transfer]\label{conj:transfer}
For $p=4/3$ and every Banach space $X$,
\[
 \calS_p(X;S^2)<\infty\quad\Longrightarrow\quad\calF_p(X)<\infty.
\]
Equivalently, for each finite $K$ there exists $C(K)<\infty$ such that
\[
 \calF_p(E)\le C(K)
 \qquad\text{for every finite-dimensional $E$ with }\calS_p(E;S^2)\le K.
\]
\end{conjecture}

The quantitative equivalence in Conjecture~\ref{conj:transfer} follows by
Gaussian $p$-moment equivalence and the exact $\ell^p$-direct-sum law for
both amplification constants: finite-dimensional counterexamples with
uniformly bounded spherical constants and divergent fractional constants
could be combined in one $\ell^p$-sum. The conjecture is a proposed
universal implication, not a consequence of the Lorentz results in this paper.

\begin{problem}[Strict-type spherical sufficiency and the reverse implication]\label{prob:reverse}
Fix any $r\in(4/3,2]$. Is
\[
 \sup_{\substack{\dim E<\infty\\T_r(E)\le K}}
 \calS_{4/3}(E;S^2)<\infty
 \quad\text{for every finite }K?
\]
More generally, does $\calF_{4/3}(X)<\infty$ imply
$\calS_{4/3}(X;S^2)<\infty$ for arbitrary Banach spaces $X$?
\end{problem}

A negative answer to the first question at one fixed $r>4/3$, with a
uniformly type-$r$ family and unbounded spherical constants, would also
answer the second question negatively, because the fractional model is
uniformly bounded on such a family by the established strict-type theorem.
On the other hand, a positive strict-type theorem would not alone settle
the reverse implication at attained type $4/3$. Neither problem is claimed
resolved here.

\subsection{Conclusion}
The Lorentz coefficient family has a sharp critical threshold and a
quantitative finite-dimensional growth law. The same logarithmic law
holds for long-memory fractional Gaussian increments and yields a sharp
spacetime bound for frequency-localized capillary waves. The remaining
open problems concern geometric structures beyond Lorentz lattices,
where Gaussian type no longer determines how curved Fourier traces
recombine.

\section*{Use of generative AI and author responsibility}
Gemini used for argument checking, and language preparation. The authors are responsible for the mathematical statements, proofs, citations, and conclusions.

\end{document}